\documentclass[12pt]{amsart}
\usepackage[english]{babel}
\usepackage[cp850]{inputenc}
\usepackage{amscd}
\usepackage{amssymb}
\usepackage{amsmath}
\usepackage{graphics}
\usepackage{graphicx}
\newtheorem{theo}{Theorem}

\newtheorem{lem}{Lemma}
\newtheorem{prop}{Proposition}
\newtheorem{coro}{Corollary}
\newtheorem{rem}{Remark}

\author{}

\begin{document}

\begin{center}
{\bf\Large  Chebyshev sets and ball operators}\\[2ex]
by\\ [2ex]
{\sc Pedro Mart\'{\i}n, Horst Martini, and Margarita Spirova}
\end{center}

\vspace*{3ex}

{\bf Abstract.}
{\small The Chebyshev set of a bounded set $K$ in a normed space is the set of centers of all minimal enclosing balls of $K$. We use the concept of ball intersection and ball hull operators to derive new properties of Chebyshev sets in normed spaces. These results give a better picture on how Chebyshev sets, ball intersections, ball hulls, and completions of bounded sets are related to each other. It is shown that the Chebyshev set of a bounded set $K$ always contains the Chebyshev set of some completion of $K$. Moreover, for a special class of sets we obtain a necessary and sufficient condition that the Chebyshev set of the respective set is a singleton. We obtain new results on critical sets of Chebyshev centers, and for that purpose, surprisingly, notions from the combinatorial geometry of convex bodies  play an essential role.   Also we give a complete geometric description of the ball hull of a  finite planar set. This can be taken  as starting point for algorithmical constructions of the ball hull of such sets.
}


\thanks{}

\thanks{Research partially supported by MICINN (Spain) and FEDER (UE) grant
MTM2008-05460, and by Junta de Extremadura grant GR10060 (partially
financed with FEDER)}
\title{}

\maketitle

\vspace{-5ex}

\section{Introduction}

\medskip

Let $K$ be a bounded set in a finite dimensional real Banach (or normed) space. The union of centers of all minimal enclosing balls of $K$ is called the Chebyshev set of $K$, and its elements are said to be the Chebyshev centers of $K$. Minimizing the maximal distance to points from $K$, Chebyshev centers occur in Approximation Theory, Banach Space Theory, Location Science, and Computational Geometry; see, e.g.,
\cite{Amir1}, \cite{Amir2}, \cite{Gar}, \cite{Gar1}, \cite{Gritz-Klee}, \cite{Al-Ma-Sp}, and \cite{Al-Ma-Sp2}.

On the other hand, the notions of ball intersection and ball hull of given point sets play an essential role in Banach Space Theory, Convex Analysis, and Classical Convexity;
see the recent contributions \cite{Mo-Sch1}, \cite{Ma-Sp}, \cite{T}, and the references given there. These concepts  are also interesting from the purely geometric point of view, because of their close relations to  notions like circumballs, minimal enclosing balls, complete sets, sets of constant width, and ball polytopes,  and also due to  their importance for constructing  such point sets; see \cite{Groe}, \cite{Al-Ma-Sp}, \cite[$\S$ 2]{MS}, \cite{Ma-Sp}, \cite{Mo-Sch2}, \cite{Mo-Sch3}, and \cite{P-M-K}.

In this paper we derive  new results on Chebyshev sets and Chebyshev centers  in real Banach spaces of finite dimension via also new  results on ball intersections and ball hulls of bounded sets in such spaces. These results give a better picture on how Chebyshev sets, ball intersections, ball hulls, and completions of bounded sets are related to each other. In particular, we show that for arbitrary sets with unique completion  the Chebyshev set and the Chebyshev radius of this set and of its completion coincide. Furthermore, we prove that the Chebyshev set of a bounded set $K$ always contains the Chebyshev set of some completion of $K$. In Section 4 we investigate the so-called critical set for a given Chebyshev center, i.e., the intersection of the boundary of the corresponding minimal enclosing ball and the original set. It turns out that this set gives  further information about the respective Chebyshev set. Surprisingly, in connection with critical sets a basic notion from the combinatorial geometry of convex bodies (namely that of  inner illuminating systems) plays an essential role. For the class of centrable sets we derive a necessary and sufficient condition that the Chebyshev set is a singleton. Finally, we give a complete geometric description of the ball hull of a finite  set in the plane. This can be used as starting point for algorithmical investigations in this direction.

\bigskip
\section{Notation and preliminaries}
\medskip

Let $\mathbb{M}^d=(\mathbb{R}^d, \|\cdot\|)$ be a $d$-dimensional real Banach space, also called a \emph{normed} (or \emph{Minkowski}) \emph{space}. The \emph{unit ball} $B$ of $\mathbb{M}^d$ is   a compact, convex set with non-empty interior (i.e., a \emph{convex body}) centered at the \emph{origin} $o$, and the boundary of $B$ is the \emph{unit sphere} of $\mathbb{M}^d$. Any  homothetical copy $x+\lambda B$ of the unit ball is called the  \emph{ball with center $x\in\mathbb{R}^d$ and radius} $\lambda > 0$ and  denoted by $B(x,\lambda)$, its boundary is the respective \emph{sphere}  $S(x,\lambda).$ We use the usual abbreviations $\mathrm{int}$, $\mathrm{relint}$, $\mathrm{bd}$, and $\mathrm{conv}$ for \emph{interior, relative interior, boundary}, and \emph{convex hull}.  The \emph{line segment} connecting  the different points $p$ and $q$ is denoted by $\overline{pq}$, the respective \emph{line} by $\langle p, q\rangle$.

\begin{figure}[h]
\begin{center}
\parbox[]{5cm}{
\scalebox{1}{\includegraphics[viewport=154 70 339 260,clip]{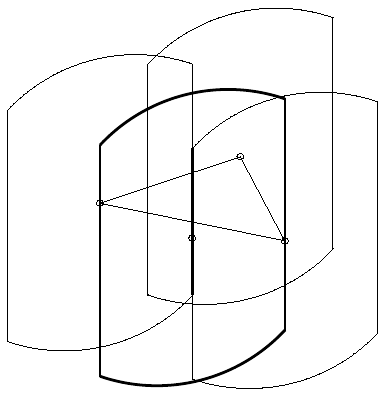}}\\}
\caption{$\mathrm{Ch }(K)$ does not belong to the convex hull of $K$}\label{fig1}
\end{center}
\end{figure}

Let $K$ be a bounded set in $\mathbb{M}^d.$ We denote  by $\operatorname{diam} (K)
:=\max\{\|x-y\|: x, y\in K\}$  the \emph{diameter of } $K$. A line segment  $\overline{pq}$ with $p, q\in K$  for which this maximum is attained is said to be a \emph{diametrical chord  of} $K$.
The \emph{$\lambda$-ball intersection} $\operatorname{bi}(K, \lambda)$ \emph{of} $K$ is the intersection of all  balls of radius $\lambda$ whose centers are in $K$:
$$\operatorname{bi}(K, \lambda)=\bigcap_{x\in K}B(x,\lambda).$$
The \emph{$\lambda$-ball hull} $\operatorname{bh}(K, \lambda)$ of $K$ is defined as the intersection of all  balls of radius $\lambda$ that contain $K$:
$$\operatorname{bh}(K, \lambda)=\bigcap_{K\subset B(x,\lambda)}B(x, \lambda).$$ Of course,   these notions make only sense if $\operatorname{bi}(K, \lambda)\neq\emptyset$ and $\operatorname{bh}(K, \lambda)\neq\emptyset$, it is clear that $\operatorname{bh}(K, \lambda)\neq\emptyset$ if and only if $\lambda\geq\lambda_K$, where $\lambda_K$ is the smallest number such that $K$ is contained in a translate of $\lambda_K B$. Such a translate is called a \emph{minimal enclosing ball of $K$} (or a \emph{circumball of $K$}), and $\lambda_K$ is said to be the \emph{minimal enclosing radius } (or  \emph{circumradius } or  \emph{Chebyshev radius}) of $K$. Clearly, we have
 \begin{equation}K_1\subseteq K_2\;\Longrightarrow \lambda_{K_1}\leq \lambda_{K_2}.\label{11} \end{equation}
 In the Euclidean subcase the minimal enclosing ball of a bounded set is always unique, but this is no longer true for an arbitrary norm. It is  easy to check that
\begin{equation}\begin{split}\{x\in \mathbb{M}^d: x\;\text{is the center of a minimal enclosing disc of}\; K\}\\=\operatorname{bi}(K, \lambda_K),\end{split}\label{19}\end{equation} yielding that
  $\operatorname{bi}(K, \lambda)\neq\emptyset$ if and only if $\lambda\geq\lambda_K$. The set of centers of minimal enclosing balls of $K$ is called the \emph{Chebyshev set of $K$}, and we denote it by $\mathrm{Ch} (K)$. Relation (\ref{19}) shows that the Chebyshev set of a bounded set is connected and convex. A recent investigation on convexity properties of Chebyshev sets in infinite dimensional normed spaces is \cite{Bor}.
   The elements of $\mathrm{Ch} (K)$ are said to be the \emph{Chebyshev centers} of $K$. Note that, in contrast to the Euclidean situation, in general the Chebyshev set of a bounded set does not  necessarily belong  to the convex hull of this set.  An example is shown in     Figure \ref{fig1}.

     According to Theorem 1 in \cite{Gar}, if every bounded set of a normed space $\mathbb{M}^d$ contains  a Chebyshev center in  its convex hull, then $d=2$ or, for $d\geq 3$,  $\mathbb{M}^d$ is Euclidean.
The best geometric description of  the Chebyshev set  of a set $K$ known until now is the fact  that $\mathrm{Ch}(K)$ is contained in  the $(\mathrm{diam} K)$-ball intersection of $K$ (see \cite[Theorem 3.3]{Bar-Pap}).

 In what follows, when we speak about the $\lambda$-ball intersection or $\lambda$-ball hull of a set $K$, we always mean that $\lambda\geq\lambda_K$.
  It is easy to check that
  \begin{equation}\label{12}\lambda_K\leq \mathrm{diam} K\leq 2 \lambda_K;\end{equation} see also \cite{Bar-Pap}.
   If $\lambda=\mathrm{diam} K$,   then we simply say  \emph{ball intersection} and  \emph{ball hull} of $K$, denoting them  by $\mathrm{bi(K)}$ and $\mathrm{bh}(K)$, respectively.

   The supremum of all $\lambda_K$, where $K$ is a bounded set of diameter 2, is \emph{Jung's constant} $J_{\|\cdot\|}$ of the considered normed space, i.e., $J_{\|\cdot\|}$ is the smallest number such that a ball of that diameter can cover any set of diameter $\leq 1$.

\smallskip

%%%%%%%%%%%%%%%%%%%%%%%%%%%%%%%%%%%%%%%%%%%%%%%%%%%%%%%%%%%%%%%%%%%%%%%%%%%%%%%%%
\bigskip
\section{Properties of the operators ball hull and  ball intersection}
\medskip

Immediately from the definition of the ball intersection we get
\begin{equation}\mathrm{bi} (\bigcup_j K_j, \lambda)=\bigcap_j \mathrm{bi}(K_j, \lambda)\label{8}\end{equation}
with $\lambda\geq\mathrm{diam} (\bigcup_j K_j)$.

Our first proposition shows that the ball intersection operator is non-increasing with respect to sets, and non-decreasing with respect to radii. Conversely, the ball hull operator is non-decreasing with respect to sets, and non-increasing with respect to radii. The Euclidean subcase of the last statement is used by Edelsbrunner (see \cite[p. 309]{Ed}) for approximating the shape of a given point set.

\begin{prop}
\begin{equation}K_1\subseteq K_2  \Longrightarrow  \operatorname{bi}(K_1, \lambda)\supseteq \operatorname{bi}(K_2, \lambda)\; and\;    \operatorname{bh}(K_1, \lambda)\subseteq \operatorname{bh}(K_2, \lambda),\label{14}\end{equation}
\begin{equation}
\lambda_1\leq \lambda_2 \Longrightarrow  \operatorname{bi}(K, \lambda_1)\subseteq \operatorname{bi}(K, \lambda_2)\; and \;   \operatorname{bh}(K, \lambda_1)\supseteq \operatorname{bh}(K, \lambda_2).\label{15}
\end{equation}
\end{prop}

\begin{proof} The inclusions in (\ref{14}) and the first inclusion in (\ref{15}) follow
directly from the definitions of ball intersections and ball hulls. Thus, only   the second inclusion in (\ref{15}) has to be proved. Let $y\in \operatorname{bh}(K,\lambda_2)$, and we  suppose that there exists $B(x,\lambda_1)$ such that $K\subset B(x,\lambda_1)$  and $y\notin B(x,\lambda_1).$ Let $x'$ be a point in $\langle x, y\rangle$ and such that $\lambda_2-\lambda_1=\|x-x'\|<\|y-x'\|.$ We prove the following for the ball $B(x',\lambda_2)$:

\begin{itemize}
\item[(i)] $B(x,\lambda_1)\subset B(x', \lambda_2)$ and

\item[(ii)]  $y\notin B(x', \lambda_2)$.

\end{itemize}
Indeed, if $z\in B(x,\lambda_1)$, then $\|z-x'\|\leq \|z-x\|+\|x-x'\|\leq\lambda_1+\lambda_2-\lambda_1=\lambda_2$  and $\|y-x'\|=\|y-x\|+\|x-x'\|=\|y-x\|+\lambda_2-\lambda_1> \lambda_2$.
Hence there exists a ball $B(x',\lambda_2)$ which contains $K$ such that $y\notin B(x',\lambda_2)$, a contradiction.\end{proof}

From \cite[Proposition 2.1.1]{Sp2} we have the following
\begin{prop}\label{1}
Let $K$ be a set in a normed space. Then
\begin{itemize}
\item[\textit{(i)}] $\operatorname{bh}(K)=\operatorname{bi}(\operatorname{bi}(K), \mathrm{diam} (K))$,
\item[\textit{(ii)}] $\operatorname{bi}(K)=\operatorname{bi}(\operatorname{bh}(K), \mathrm{diam} (K))$.
\end{itemize} \end{prop}

Note that if $\lambda\geq\mathrm{diam} K$, then we have the inclusions
\begin{equation} K\subseteq \operatorname{bh}(K, \lambda) \subseteq \operatorname{bi}(K, \lambda).\label{6}\end{equation}
The first inclusion follows immediately from the definition of the ball hull. For the second inclusion we take $x\in\mathrm{bh}(K, \lambda)$. Thus $x$ belongs to all balls $B(y, \lambda)$ which contain $K$. If $y$ is an arbitrary point from $K$, then $B(y, \lambda)\supseteq B(y, \mathrm{diam} (K))\supseteq K$. Therefore every ball with radius $\lambda$ and center from $K$ contains $K$.

If $\lambda=\mathrm{diam} K$, then
\begin{equation} K\subseteq \operatorname{bh}(K)\subseteq K^c \subseteq \operatorname{bi}(K),\label{7}\end{equation}
where $K^c$ is a \emph{completion} of $K$, i.e., a complete set that contains $K$ and has the same diameter. (Recall that a set in $\mathbb{M}^d$ is said to be \emph{complete} if it cannot be enlarged by additional points  without increasing its diameter.) The inclusion $\operatorname{bh}(K, \lambda)\subseteq K^c$ was given in the proof of Theorem 5 in \cite{Groe}.
The completeness of a set $K$ is equivalent of the fact that $K=\mathrm{bi}(K)$; see \cite[p. 167, (E)]{E2}.   It   also known  that  completeness and constant width are equivalent in the Euclidean subcase and in every two-dimensional normed space. (For the definition and many properties of \emph{sets of constant width} in normed spaces we refer to Section 2 of the survey \cite{MS}.) In a normed space of dimension at least three, constant width implies completeness, but not vice versa; cf. \cite{E2} and  \cite[pp. 98-99]{MS}.

The next statement shows relations  between the notions of ball hull, ball intersection and  completion of a given set.

\begin{prop}\label{theo3a}
Let $K$ be a bounded set   in a normed space. Then
\begin{equation}\label{9}
\mathrm{bh}(K)=\bigcap(\text{all completions of} \;K)
\end{equation}
and
\begin{equation}\label{10}
\mathrm{bi}(K)=\bigcup(\text{all completions of} \;K).
\end{equation}
\end{prop}

\begin{proof}
First we prove (\ref{10}). Due to (\ref{7})
we only need to check whether  for every point $x$ of $\mathrm{bi}(K)$ there exists a completion of $K$ such that  $x$ belongs to this completion.
%Let us assume that for all completions $K^c$ of $K$ the point $x$ does not belong to $K^c$.
Since $x\in \mathrm{bi}(K)$, for any point $y$ of $K$ the inequality $\|x-y\|\leq \mathrm{diam}(K)$ holds. This means that $\mathrm{diam}(K\bigcup \{x\})=\mathrm{diam}(K)$, and there exits a completion $K^c$ of the set $K\bigcup\{x\}$ which is also a completion of $K$.

Statement (\ref{9}) follows from (\ref{10}) and (\ref{8}).
\end{proof}

\begin{rem}\label{rem2} The Euclidean planar subcase  of Proposition \ref{theo3a}  was proved  by Bavaud; see \cite[Theorem IV]{Bav}. Another proof of this statement was given by Moreno in \cite{Mo}; see  Proposition 2 there. Note that in that paper,  as well as in the  papers of Moreno and  Schneider cited here, the notation $\eta(K)$ and $\theta(K)$, respectively, for  ball intersection and  ball hull of $K$ is used. In fact,
 $\theta(K)$ is defined as $\cap_{x\in \eta(K)}B(x, \mathrm{diam} K)$. Due to Proposition  \ref{1}, (i), the latter is  the ball hull $\mathrm{bh}(K)$. In \cite{Mo-Sch2} and \cite{Mo-Sch3} Moreno and Schneider  call $\mathrm{bi}(K)$ and $\mathrm{bh}(K)$, respectively, the \emph{wide spherical hull} and the \emph{tight spherical hull} of $K$. Also we note  that Baronti and Papini used in \cite{Bar-Pap} the notation $K'$ for ball intersection  and $K^c$ for  ball hull of $K$.
\end{rem}

The first corollary of  Proposition \ref{theo3a} follows immediately from (\ref{10}).

\begin{coro}\label{cor1}
Every ball with radius $\lambda=\mathrm{diam} K$ and center from $K$ which contains $K$  also contains any completion of $K$.
\end{coro}

As it was announced in Section 2, we have that $\mathrm{Ch}(K)\subset\mathrm{bi}(K)$; see \cite[Theorem 3.3]{Bar-Pap}. This fact and  (\ref{10}) from Proposition  \ref{theo3a} yield
\begin{coro}\label{cor2}
Let $B(x, \lambda_K)$ be a minimal enclosing ball of a set $K$. Then there exists a completion  $K^c$ of $K$ which contains $x$.
\end{coro}

For the next considerations we need  Theorem 3.5 from \cite{Bar-Pap} which says: for every ball  of radius $\lambda\geq\frac{1}{2}\mathrm{diam} K$ which contains $K$
there exists a completion $K^c$ of $K$ contained in this ball. This theorem immediately implies

\begin{lem}\label{theo4}
Let $B(x, \lambda_K)$ be a minimal enclosing ball of $K$. Then $B(x, \lambda_K)$ contains some completion $K^c$ of $K$.
\end{lem}

%\begin{theo} Let $B(y, \lambda_K)$ be a minimal enclosing ball of $K$. Then $B(y, \lambda_K)$ contains every completion $K^c$ of $K$.
%\label{theo4}
%\end{theo}

%\begin{proof}
%Assume that there exists a point $x$ from some completion $K_c$ with $x\not\in B(y, \lambda_K)$. Let the ray  emanating from $y$  and passing through $x$ intersect the boundary of $B(y, \lambda_K)$ in $p$, and $H_p$ be a supporting hyperplane of $B(y, \lambda_K)$ at $p$. Let $\varphi$ be the homothety with center $p$ mapping $B(y, \lambda_K)$ into $B(\varphi(y), \mathrm{diam} K)$. Since $\mathrm{diam} K\geq \lambda_K$, we have
%$B(\varphi(y), \mathrm{diam} K)\supseteq B(y, \lambda_K)\supset K$, which implies that $\mathrm{bh} (K)\subset B(\varphi(y), \mathrm{diam} K)$. Since every completion of $K$ is contained in its ball hull (Theorem \ref{theo3}, (\ref{9})), we obtain  the inclusion  $K^c\subset B(\varphi(y), \mathrm{diam} K)$ which contradicts the fact that $H_p$ separates $B(\varphi(y), \mathrm{diam} K)$ and $x\in K^c$.
%\end{proof}

\begin{prop}\label{cor3}
Every minimal enclosing ball of a set $K$ contains its ball hull.
\end{prop}

\begin{proof}
This  follows from Lemma \ref{theo4} and (\ref{7}).
\end{proof}

\begin{lem}\label{cor4}
For any bounded set $K\subset\mathbb{M}^d$ there exists some completion $K^c$ such that the relations $\lambda_K=\lambda_{K^c}$ and  $\mathrm{Ch }(K^c)\subseteq\mathrm{Ch }(K)$ hold.
\end{lem}

\begin{proof}
By Lemma \ref{theo4} and inequality (\ref{11}) there exists a completion $K^c$ with the first relation. In view of  $\mathrm{Ch }(K)=\mathrm{bi}(K, \lambda_K)$,  $\mathrm{Ch }(K^c)=\mathrm{bi}(K^c, \lambda_K)$, and (\ref{14}) we get
$\mathrm{Ch }(K^c)\subseteq\mathrm{Ch }(K)$.
\end{proof}

\begin{lem}\label{cor4.1}
For any bounded set $K\subset\mathbb{M}^d$ and $x\in \mathrm{Ch }(K)$ there exists some completion $K^c$ such that $x\in\mathrm{Ch }(K^c)$.
\end{lem}

\begin{proof}
Let $x\in \mathrm{Ch }(K)$. By Lemma \ref{theo4} there exists a completion $K^c$ contained in $B(x,\lambda_K)$ with $\lambda_{K^c}=\lambda_K.$ Therefore, $\|x-y\|\leq \lambda_{K^c}$ for every $y\in K^c$ and $x\in \cap_{y\in K^c} B(y,\lambda_{K^c})$. We conclude that
$x\in \mathrm{bi }(K^c, \lambda_{K^c})=\mathrm{Ch }(K^c)$.
\end{proof}

The next theorem gives a relation between the Chebyshev set of a set and the Chebyshev sets of its completions.
\begin{theo}\label{cor5}
For any bounded  set $K\subset\mathbb{M}^d$  the inclusions
$$\bigcap \mathrm{Ch }(K^c) \subseteq \mathrm{Ch }(K)\subseteq \bigcup \mathrm{Ch }(K^c)$$ hold.
\end{theo}

\begin{proof}
The first inclusion follows by Lemma \ref{cor4} and the second one by Lemma \ref{cor4.1}.

\end{proof}
\begin{figure}[ht]
\begin{center}
\hspace*{-0.5cm}
\parbox[]{5cm}{
\includegraphics[width=4.5cm]{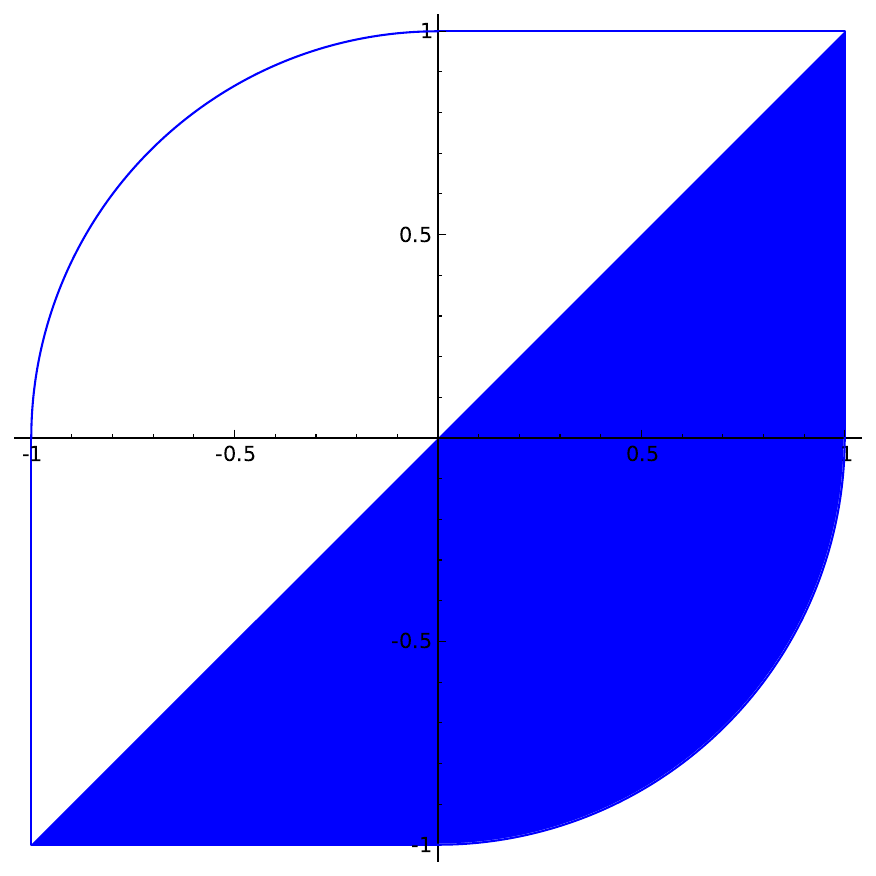}\\}
\hspace{0.5cm}
\parbox[]{6cm}{
\scalebox{0.65}{\includegraphics[viewport=6 145 271 416,clip]{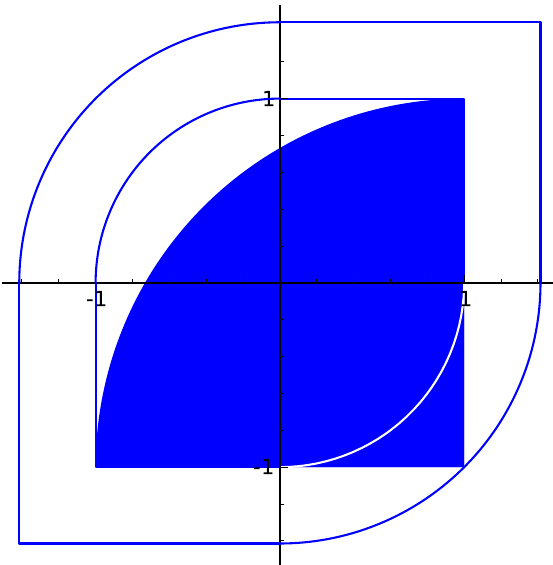}}\\}
\caption{A set $K$ and a completion $K^c$}\label{fig2.1}
\end{center}
\end{figure}
Note that the equality $\mathrm{Ch }(K)= \cup \mathrm{Ch }(K^c)$ is not true in general. If we consider in $\mathbb{R}^2$ the normed space with unit ball as in Figure \ref{fig2.1} (left), then the Chebyshev  set of the set $K=\{(x,y)\in B(0,1)\, : x-y<0\}$ is the point $(0,0)$, but the Chebyshev set of the completion $K^c$ shown in Figure \ref{fig2.1} (right) is not  this point. It should also be  noticed that the location of the Chebyshev set of a bounded set as given in Theorem \ref{cor5} sharpens  the inclusion $\mathrm{Ch}(K)\subset \mathrm{bi}(K)$ given in  \cite[Theorem 3.3]{Bar-Pap}. Indeed, by  $\mathrm{Ch}(K)\subset \mathrm{bi}(K)$ and (\ref{14}), we have
\begin{align*}
\bigcup \mathrm{Ch}(K^c) & \subseteq  \bigcup\mathrm{bi}(K^c, \mathrm{diam}(K^c))\\
& =\bigcup\mathrm{bi}(K^c,\mathrm{diam}(K))\subseteq\bigcup \mathrm{bi}(K)=\mathrm{bi}(K).
\end{align*}

The last corollary of this section follows immediately from  Theorem \ref{cor5}.

\begin{coro}\label{cor10} If $K$ is a bounded set with unique completion $K^c$, then the relations $\lambda_K=\lambda_{K^c}$ and $\mathrm{Ch}(K)=\mathrm{Ch}(K^c)$ hold.
\end{coro}

It should be noticed that the class of sets with unique completion is relatively large  and well studied; see \cite{MS} and the references given there.
We hope that Corollary \ref{cor10} will be of importance for  further investigations  of Chebyshev sets and  Chebyshev radii of sets from this class. Due to this corollary, one does not have to investigate these notions for an arbitrary set from this class. It is sufficient to study them only for complete sets.

\begin{figure}[h]
\begin{center}
\includegraphics[width=6cm]{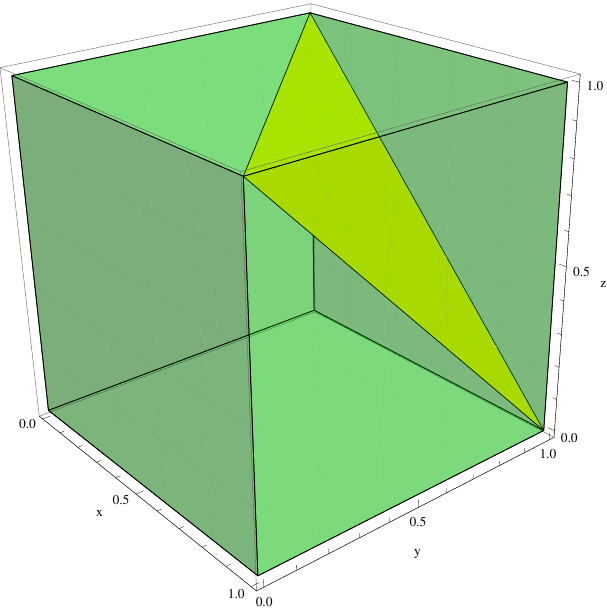}\\
\caption{$\mathrm{Ch }(K)$ does not intersect the convex hull of $K$}\label{fig2}
\end{center}
\end{figure}

\begin{rem}
We would like to comment  the chain  (\ref{7}) of inclusions as follows. For a bounded convex set $K$, in general $Ch(K)\nsubseteq K$; see again Theorem 1 in \cite{Gar}. According to \cite[Theorem 3.3]{Bar-Pap} we have that $\mathrm{Ch}(K)\subseteq\mathrm{bi}(K)$ and, more precisely, Corollary \ref{cor2} implies that there exists a completion $K^c$ of $K$ with $\mathrm{Ch}(K)\subset K^c$. The following example shows that $\mathrm{Ch}(K)\nsubseteq\mathrm{bh}(K)$. Moreover, in  general normed spaces there exist sets  whose Chebyshev sets do not  necessarily intersect their convex hull. In Figure \ref{fig2}  a set with this property is shown. Consider the maximum  norm in $\mathbb{R}^3$. Then  the set $K=\{(1,1,0), (1,0,1), (0,1,1)\}$ has only one Chebyshev center which does not belong  to the convex hull of these points. It is easy to construct  a similar example with non-empty interior.

\end{rem}

%\begin{coro}\label{cor5}
%The Chebyshev set of a set $K$ belongs to the ball hull of $K$.
%\end{coro}

%\begin{proof}

%If we prove  that every completion $K^c$ of $K$ contains $\mathrm{Ch} (K)$ then, by Proposition \ref{theo3a}, (\ref{9}), the proof is complete.
%On the contrary, we assume that there exists a Chebyshev center $x$ not being  from $K^c$. According to Corollary \ref{cor4} the ball $B(x, \lambda_K)$ is also a minimal enclosing ball of $K^c$ which implies that $K^c\cap \mathrm{bd} B(x, \lambda_K)\neq \emptyset$. If $y\in K^c\cap \mathrm{bd} B(x, \lambda_K)$, then $\lambda_K=\|x-y\|> \mathrm{diam} K$, a contradiction to  (\ref{12}).
%\end{proof}

%%%%%%%%%%%%%%%%%%%%%%%%%%%%%%%%%%%%%%%%%%%%%%%%%%%%%%%%%%%%%%%%%%%%%%%%%%%%%%%%%%%%%
\bigskip
\section{Critical sets}
\medskip

Let $B(x, \lambda_K)$ be a minimal enclosing ball of a set $K$. The set of all boundary points of $K$ which belong to $\mathrm{bd} B(x, \lambda_K)$ is called the \emph{critical set} of $K$ with respect to $x$ and denoted by $\mathrm{cr}(K, x)$. The elements of $\mathrm{cr}(K, x)$ are said to be \emph{critical points}. The notion of critical set was introduced by Garkavi in \cite{Gar1}, and it turns out that  such sets yield  important information on Chebyshev sets.  E.g., via critical points   we can completely clarify, for the class of centrable sets, whether or not  the Chebyshev set is a singleton; see Proposition \ref{pro3}.  Sets $K$ with the property $\mathrm{diam} K=2\lambda_K$ are called \emph{centrable sets}.  In other words,  a set is centrable if   the second inequality in (\ref{12}) is an equality for it. Centrable  sets are treated, e.g., in \cite{Rao}, where  real Banach spaces with the property that every finite set is centrable are characterized.

\begin{theo}\label{theo5} Let there be given a non-centrable compact convex set $K$ in a normed plane $\mathbb{M}^2$, and $B(x,\lambda_K)$ be a minimal enclosing ball of $K$.
 Then there exist $m\in\{2, 3\}$ points $t_1,\dots,t_m\in \operatorname{cr}(K, x)$ being the vertices of an $(m-1)$-simplex such that
$$
x\in \operatorname{relint}(\operatorname{conv}(t_1,\dots,t_m)).
$$
\end{theo}

\begin{proof}
By Lemma 2.2 in \cite{Br-Koe}, there exist  $m$ critical points $t_1,\dots, t_m$ ($m\in \{2,3\}$) such that the minimal enclosing radius of the set $S=\operatorname{conv}(t_1,\dots,t_m)$ is $\lambda_K.$

\begin{figure}[ht]
\begin{center}
\scalebox{1}{\includegraphics[viewport=63 404 198 564,clip]{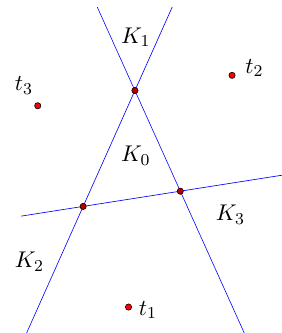}}\\
\caption{}\label{fig3}
\end{center}
\end{figure}

If $m=2$ and $S=\operatorname{conv}(t_1,t_2)$, then $x\in B(t_1,\lambda_K)\cap B(t_2,\lambda_k)$. The intersection of these two balls has no interior points, because in that case there exists a ball  containing $t_1$ and $t_2$ with radius smaller than $\lambda_K$. Therefore (see \cite[Section 3.3]{MSW}), the intersection is a segment and $K$ is centrable due to $\|t_1-t_2\|=2\lambda_K$.

\begin{figure}[h]
\begin{center}
\scalebox{1}{\includegraphics[viewport=87 509 207 697,clip]{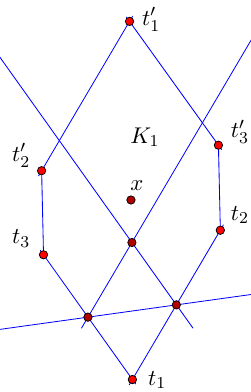}}\\
\caption{}\label{fig4}
\end{center}
\end{figure}
Let us assume that $m=3$ and $S=\operatorname{conv}(t_1,t_2,t_3)$. The center $x$ of any  minimal enclosing ball can belong to some of the regions $K_0$, $K_1$, $K_2$ or $K_3$ shown in  Figure \ref{fig3}; these are bounded by the lines meeting the midpoints of $\overline{t_1,t_2}$, $\overline{t_1, t_3}$ and $\overline{t_2, t_3}$ (see \cite{Al-Ma-Sp} and \cite{Al-Ma-Sp2}). If $x\in  K_0$, the proof is finished. Let us assume that $x$ belongs to  $K_1$ (the cases $x\in K_2$ or $x\in K_3$ are similar to this one). The situation is like in Figure \ref{fig4} (see again \cite{Al-Ma-Sp} and \cite{Al-Ma-Sp2}), where $t_i'$ ($i\in \{1,2,3\}$) is  symmetric to $t_i$ with respect to $x$. Therefore  $\operatorname{conv}(t_1,t_2,t_3,t_1',t_2',t_3')$ belongs to $B(x,\lambda_K$).

 Let $x'$ be the translate of  $x$  in  direction  $t_3-t_2'$  such that  $x'$ lies on  $\langle t_2, t_3\rangle$. It is clear that  $\|t_1-t_2\|\leq 2\lambda_K$. $K$ would be  centrable if   equality holds. In the other case,  $B(x',\lambda_K)$ contains $S$ and there exists a ball of radius $\lambda<\lambda_K$ which also contains $S$. Therefore, $\lambda_K$ is not the minimal enclosing radius of $S$, a contradiction.

\end{proof}

As an immediate consequence of Theorem \ref{theo5} we have a statement which  completely clarifies  the location of Chebyshev sets of non-centrable planar sets.

\begin{coro}
\label{cor9} If $K$ is a non-centrable set in a normed plane,  then the Chebyshev set of $K$ belongs to the relative interior of the convex hull of $K$.
\end{coro}

The following notion from the combinatorial geometry of convex bodies is also related to our investigations:

A subset $P$ of the boundary of a convex body $K$  is called  an \emph{inner illuminating system} of $K$ if for every $x\in \mathrm{bd} K$ there exists at least one $p\in P$ such that the line segment $\overline{px}$ meets the interior of $K$; cf. \cite[$\S$ 34]{Bo-Ma-So}.

\begin{figure}[h]
\begin{center}
\scalebox{.7}{\includegraphics[viewport=111 491 365 754,clip]{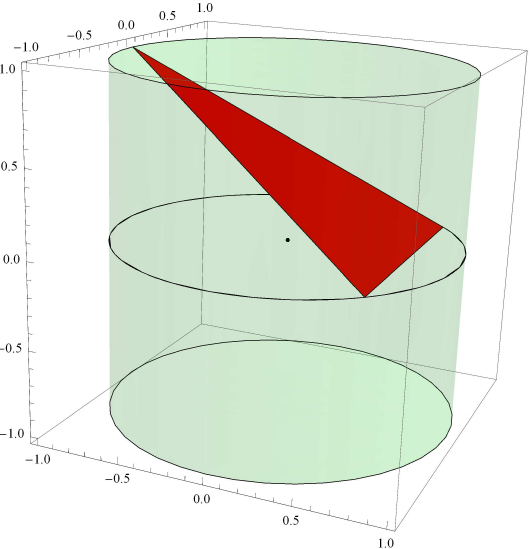}}\\
\caption{$\mathrm{Ch }(K)$ does not belong to $\mathrm{relint}(K)$}\label{c3}
\end{center}
\end{figure}
  As Figure \ref{c3} shows,  Theorem \ref{theo5} is not true in dimensions at least three. Moreover, it is not true that every convex body has a minimal enclosing ball with the property described in Theorem \ref{theo5}. Otherwise, every bounded set would have a Chebyshev center belonging to its convex hull. This contradicts  Theorem 1 in \cite{Gar} saying that this is only possible when the underlying normed space is of dimension 2 or, for any dimension $d\geq 3$, when it is Euclidean. But if for a Chebyshev center  points $t_1, \ldots, t_m$ exist as in Theorem \ref{theo5},  then these points form an inner illuminating system of the respective convex body.  Since such  systems are well studied in the combinatorial geometry of convex bodies, we are motivated for the next definition. Let $K$ be a convex body in a normed space $\mathbb{M}^d$  with Chebyshev center $x$ having the following  property: for some $m\in\{2, \ldots, d+1\}$ there exist $m$ points $t_1, \ldots, t_m\in \mathrm{cr}(K, x)$ being the vertices of an $(m-1)$-simplex such that
$$x\in \operatorname{relint}(\operatorname{conv}\{t_1,\dots,t_m\}).$$ Then $x$ is called a \emph{meanstream  Chebyshev center} of $K$. Our Theorem \ref{theo5} shows that all Chebyshev centers in the planar case are meanstream. If $x$ is a meanstream Chebyshev center of $K$, then the smallest integer $m$  with the above property is   called    the \emph{degree} of $x$, and the system $\{t_1, \ldots, t_{m}\}$ is said to be a \emph{base system} of $x$.

The next proposition follows directly from the definition of the degree  of a meanstream Chebyshev center.

\begin{prop}\label{pro1}
Let $K$ be a bounded  set in a normed space $\mathbb{M}^d$. Then $K$ is centrable  if and only if it has a meanstream Chebyshev center of degree 2.
\end{prop}

 If $K$ is  centrally symmetric, then a subset of $\mathrm{bd} K$  is said to be \emph{global} if it is not contained in any closed halfspace whose bounding hyperplane passes through the center of $K$; see \cite{Grue}.

\begin{prop}\label{cor6} Let $K\subset (\mathbb{R}^d, \|\cdot\|)$ be a convex body having  a meanstream Chebyshev center $x$  of degree $m$ with  base system $\{p_1, \ldots, p_{m}\}$. Then
the set  $\{p_1, \ldots, p_{m}\}$ is an inner illuminating system of the minimal enclosing ball $B(x, \lambda_K)$ of $K$.
\end{prop}

\begin{proof} We fix some point $p_i$. The only points of $\mathrm{bd} B(x, \lambda_K)$ not  inner illuminated by $p_i$ are those which belong to a proper face of $B(x, \lambda_K)$. But not all $p_1, \ldots, p_{m}$ belong to such a proper face.
\end{proof}

Analogously we get

\begin{prop}\label{cor7} Let $K\subset (\mathbb{R}^d, \|\cdot\|)$ be a convex body with meanstream  Chebyshev center $x$ of degree $d+1$ and  base system $\{p_1, \ldots, p_{d+1}\}$.  Then $\mathrm{cr}(K, x)$ is global, and  $\{p_1, \ldots, p_{d+1}\}$ is an inner illuminating system  of $K$.
\end{prop}

% Sets $K$ with property $\mathrm{diam} K=2\lambda_K$ are called \emph{centrable sets}.  In other words,  a set is centrable if   the second inequality in (\ref{12}) is an equality for it. Centrable  sets are treated, e.g., in %\cite{Rao}, where  real Banach spaces with the property that every finite set of such a space  is centrable are characterized.
 It is now our aim to clarify  when, for a given centrable set $K$,  $\mathrm{Ch}(K)$ is a singleton. For that reason we need some basic notions referring to the face structure of a convex body. First we recall that a boundary point $q$ of a convex body $K$ is called an \emph{extreme point} if it is not in the relative interior of any line segment belonging to  $K$.
%If $q\in \mathrm{bd} K$ is not extreme, then the set of all segments from the boundary of $K$ having $q$ in their relative interior is said to %be a \emph{face} of $K$.
% It is clear that if $q$ is an extreme point of $K$ and there exists a face of $K$ containing $q$, then $q$ does not belong to the relative %interior of this face.
%Note also that if $q$ is an extreme point does not belonging to any face of $K$, then there exists no line segment of the boundary of $K$ %containing $q$. Thus we have

It should be noticed that if $x$ is a meanstream  Chebyshev center of a centrable set $K$ with base $\{p_1, p_2\}$, then $x=\frac{1}{2}(p_1+p_2)$ and $\lambda_K=\frac{1}{2}\|p_1-p_2\|$.

\begin{prop}\label{pro2}
 The Chebyshev set of two different points $p_1$ and $p_2$  is a singleton if and only if $p_1$ (and therefore $p_2$) is an extreme point of $B(\frac{1}{2}(p_1+p_2), \frac{1}{2}\|p_1-p_2\|)$.
\end{prop}

\begin{proof}
Let $p_1$ be an extreme point. Assume that there exists a Chebyshev center $x+v$, $v\neq 0$, of $\{p_1, p_2\}$ different to  $x$. Then
$$2\lambda_K=\|p_1-p_2\|\leq \|p_1-(x+v)\|+\|p_2-(x+v)\|\leq 2\lambda_K,$$ where $K=\{p_1, p_2\}$.
This is only possible when $\|p_1-(x+v)\|=\|p_2-(x+v)\|=\lambda_K$. Hence the points $p_1-v$ and $p_1+v$ are from  $\mathrm{bd} B(x, \lambda_K)$, a contradiction. The  proof of the converse  can be  omitted.

\end{proof}

\begin{prop}\label{pro3}
 Let $K$ be a centrable set in a normed space $\mathbb{M}^d$. The Chebyshev set of $K$ is a singleton if and only if  $K$ has a Chebyshev center $x$ with base  system $\{p_1, p_2\}$ such that $p_1$ (and therefore $p_2$) is an extreme point of $B(x, \lambda_K)$.
\end{prop}

\begin{proof}
%If $K^c$ is a completion of $K$,  then it is a completion of $\{p_1, p_2\}$, too. Thus Corollary \ref{cor4} yields
%$\mathrm{Ch}\{p_1, p_2\}=\mathrm{Ch}(K)=\mathrm{Ch}(K^c)$.
Since $K$ is centrable, by Proposition \ref{pro1} there exist a Chebyshev center $x$ and two points $\{p_1,p_2\}\in \mathrm{cr} (K,x)$ such that $x\in \operatorname{relint}(\operatorname{conv}\{p_1,p_2\}).$ Therefore, $x=\frac{p_1+p_2}{2}$, $\lambda_K=\frac{\|p_1-p_2\|}{2}$ and $p_1$ is the point symmetric to  $p_2$ with respect to $x$.

Let us assume that $\mathrm{Ch }(K)$ is a singleton. If $p_1$ and $p_2$ are not extreme points of $B(x,\lambda_K)$, then both $p_1$ and $p_2$ belong to parallel segments with direction $v$ which are contained in the boundary of $B(x,\lambda_K)$. Let us consider  parallel support hyperplanes of $B(x,\lambda_K)$ at $p_1$ and at $p_2$ which contain the direction $v$. Moving the center $x$ in the direction $v$, it is possible to find a new Chebyshev center of $K$, a contradiction.

Let us assume that $p_1$ and $p_2$ are extreme points of $B(x,\lambda_K)$. Since $x=\frac{p_1+p_2}{2}$ and $\lambda_K=\frac{\|p_1-p_2\|}{2}$, Proposition \ref{pro2} yields $\mathrm{Ch }(\{p_1,p_2\})=x$ and $\lambda_{\{p_1,p_2\}}=\lambda_K$. As $\{p_1,p_2\}\subseteq K$, by (\ref{14}) we obtain that
$$x=\mathrm{Ch }(\{p_1,p_2\})=\mathrm{bi }(\{p_1,p_2\},\lambda_K)\supseteqq \mathrm{bi }(K,\lambda_K)=\mathrm{Ch }(K)$$
\end{proof}

%%%%%%%%%%%%%%%%%%%%%%%%%%%%%%%%%%%%%%%%%%%%%%%%%%%%%%%%%%%%%%%%%%%%%
%\section{Sets whose ball intersection or ball hull are compete}

%\begin{theo}\label{theo2}
%Let $K$ be a convex body  in a normed  space $\mathbb{M}^d$. Then  the following conditions are equivalent:
%\begin{enumerate}
%\item[(i)] $\operatorname{diam}(\operatorname{bi}(K))=\mathrm{diam} K$,
%\item[(ii)] $\operatorname{bi}(K)=\operatorname{bh}(K)$,
%\item[(iii)] $\operatorname{bi}(K)$ is complete,
%\item[(iv)] $\operatorname{bh}(K)$ is complete,
%\item[(v)] $K$ has a unique completion.
%\end{enumerate}
%\end{theo}

%\begin{proof} See \cite[Theorem 5]{Groe} and \cite[Theorem 3.7]{Bar-Pap}.

%\end{proof}
%\smallskip

%%%%%%%%%%%%%%%%%%%%%%%%%%%%%%%%%%%%%%%%%%%%%%%%%%%%%%%%%%%%%%%%%%%%%%%%%
\bigskip

\bigskip
\section{On ball hulls of  finite  sets in the planar case}
\medskip

As  Proposition \ref{cor3} shows, the ball hull of a set occurs when studying minimal enclosing balls, giving  also motivations for  more precise descriptions of ball hulls of  planar finite sets. Some authors make use of the ball hull in the Euclidean plane (also known as the $\alpha$-hull or the circular hull) for developing algorithms in combinatorial geometry (e.g., \cite[p.309]{Ed}) and cluster analysis (e.g., \cite{H}, \cite{H2}, \cite{Sh}). In this section, we hope to present the foundations of similar algorithms for general norms. For our considerations we need the   following  lemma which  goes back to Gr\"{u}nbaum \cite{Grue1} and Banasiak \cite{Ban}; see also \cite[$\S$ 3.3]{MSW}.

\bigskip

\begin{lem}\label{twocircles} Let $\mathbb{M}^2$ be a normed plane. Let $C\subset \mathbb{M}^2$ be a compact, convex disc whose  boundary is the closed curve $\gamma$; $v$ be a vector in $ \mathbb{M}^2$; $C+v$ be a translate of $C$ with boundary $\gamma'$. Then $\gamma\cap \gamma'$ is the  union of two segments, each of which may degenerate to a point or to the empty set.

Suppose that this intersection consists of two connected non-empty components $A_1$, $A_2$. Then the two lines of translation supporting $C\cap C'$ intersect $C\cap C'$  exactly in $A_1$ and $A_2.$

Choose a point $p_i$ from each component $A_i$ and let $c_i=p_i-v$ and $c_i'=p_i+v$ for $i=1,2.$ Let $\gamma_1$ be the part of $\gamma$ on the same side of the line $\langle p_1, p_2\rangle$ as $c_1$ and $c_2$; let $\gamma_2$ be the part of $\gamma$ on the  side of  $\langle p_1, p_2\rangle$ opposite to  $c_1$ and $c_2$; and similarly for $\gamma$, $\gamma_1'$, and $\gamma_2'$.

Then $\gamma_2\subseteq \operatorname{conv}(\gamma_1')$ and $\gamma_2'\subseteq \operatorname{conv}(\gamma_1).$
\end{lem}

Let $p$ and $q$  be two points of  the circle $S(x,\lambda)$. In the following, the \textit{minimal circular arc of $B(x,\lambda)$ meeting $p$ and $q$ } is the piece of $S(x,\lambda)$  with endpoints $p$ and $q$ which lies in the half-plane bounded by the line $\langle p, q\rangle$ and   not containing the center $x$. If $p$ and $q$ are opposite in $S(x,\lambda)$, then the two half-circles with endpoints $p$ and $q$ are minimal circular arcs of $S(x,\lambda)$ meeting $p$ and $q$. In the following, we suppose for simplicity $\lambda=1$. We denote a minimal circular arc meeting $p$ and $q$ by  $\widehat{pq}$.

\smallskip

With Lemma \ref{twocircles} we can prove a  generalization of  Lemma 4.1 in \cite{H}.

\begin{lem}\label{3}
Let $\mathbb{M}^2$ be a normed plane with unit disc $B$ and  points $p, q\in B$ satisfying $\|p-q\|\leq 1$.
 \begin{enumerate}
 \item If  $p, q\in S(o,1)$ and there exists another circle $S(x,1)$ through $p$ and $q$, then either  the origin $o$ and $x$ are in different half-planes with respect to the line $\langle p, q\rangle$, or the segment $\overline{pq}$ belongs to $S(x,1)\cap S(o,1)$.

 \item Any minimal circular arc of radius 1 meeting $p$ and $q$ also belongs to $B$.

 \item  If a circular arc meeting $p$ and $q$ is  contained in $B$ such that it contains interior points of $B$, then this arc is a minimal circular arc.

 \end{enumerate}
\end{lem}

\smallskip

\begin{proof} $(1)$ Let $v$ be the vector $p-q$. Let us assume that $x$ and $o$ are in the same half-plane  with respect to the bounding  line $p+\lambda v$ (if $x$ and $o$ are in different half-planes defined by the line $p+\lambda v$, $\lambda\in \mathbb{R}$, then nothing has  to  be proved).

Case 1: The point \emph{$x$ lies between the lines $o+\lambda v$ and $p+\lambda v$}, $\lambda\in \mathbb{R}$.  The point $x$ is not from $\overline{oq}$ since $\|q-o\|=\|q-x\|=1$. For  similar reasons, $x$ is not from  $\overline{op}$. If $x$ is in the interior of the triangle defined by $o$, $p$ and $q$, then $S(p,1)$ is not a convex curve, since $o$ and $x$ are in $S(p,1)$ and $q\in B(p,1)$.  If $x$ is not in the triangle defined by $o$, $p$ and $q$, let us suppose that $x$ and $p$ are in different half-planes  with respect to  $\langle o, q\rangle$ (the argument is similar if $x$ and $q$ are in different half-planes  with respect to  $\langle o, p\rangle$). The points $x+p$, $x+q$, $p$, and $q$ are in $S(x,1$), but the two  pairs of points $(x+p, x+q)$ and $(p, q)$ define two parallel segments. Therefore the four points are in the same line, and the segment $\overline{pq}$ is in $S(x,1)$.

Case 2: The origin \emph{$o$ lies between the lines $x+\lambda v$ and $p+\lambda v$, $\lambda\in \mathbb{R}$}. We can argue like in Case 1, exchanging the role of $x$ and $o$, because the roles that $o$ and $x$  play in this part of the Lemma are symmetric.

 $(2)$ Let us consider a minimal circular arc of radius 1 meeting $p$ and $q$, and let $S(x,1)$ be the circle that contains this arc. The curves $C=S(o,1)$ and $C'=S(x,1)$ verify the hypothesis of Lemma \ref{twocircles}. Let $p_i$, $\gamma_i$ and $\gamma_i'$ (with $i=1,2$) be as in Lemma \ref{twocircles}. There exist two points $p_1$ and $p_2$ such that the component $\gamma_2'$ is maximal. Then the points $p$ and $q$ and the minimal circular arc meeting $p$ and $q$ are in $\gamma_2'$. By Lemma \ref{twocircles} we have $\gamma_2'\subseteq \operatorname{conv}(\gamma_1)\subseteq B(0,1).$

$(3)$ Let us assume that $p$ and $q$ are in $B$, $\|p-q\|\leq 1$, and that there exists an arc meeting $p$ and $q$  which has points in  $\mathrm{int} B$. Let $B(x,1)$ be the ball such that the arc is in $S(x,1)$.  We have $x\neq o$, because there are points of the arc in  $\mathrm{int} B$. The set  $S(o,1)\cap S(x,1)$ is the union of two segments, each of which may degenerate to a point or to the empty set. If the intersection consists of two connected components $A_1$ and $A_2$, we can choose two points $p_1\in A_1$ and $p_2\in A_2$ like in Lemma \ref{twocircles} such that there is a maximal component $\gamma_2$ of $S(x,1)$ defined by $p_1$ and $p_2$ and   satisfying the following: it is inside $B$, and the other component of $S(x,1)$ has only the two points $p_1$ and $p_2$ in $B$. The points $p$ and $q$ are in $B$,  and the minimal arc of $S(x,1)$ meeting $p$ and $q$ is also in $B$,  by the above result. Then this minimal arc is in $\gamma_2$. Therefore, the other arc defined in $S(x,1)$ by $p$ and $q$ has points outside of  $B$.

\end{proof}

Let $p$ and $q$ be two points such that $\|p-q\|\leq 1$. By  part $(1)$ of Lemma \ref{3}, there exist only two minimal arcs meeting $p$ and $q$ (which may degenerate to the segment $\overline{pq}$ if the centers of the unit balls are in the same half-plane).

From part $(2)$ of Lemma \ref{3} we get immediately

\bigskip

\begin{lem}\label{4}
Let $K=\{p_1,p_2,\dots,p_n\}$ be a finite set in a normed plane   $\mathbb{M}^2$ having  diameter 1. Then any possible unit ball $B(x,1)$ which contains $K$ also contains every minimal circular arc of radius 1 meeting $p_i$ and $p_j$ ($i,j=1,2,\dots, n$).
\end{lem}

We continue with a further lemma.

\bigskip

\begin{lem}\label{lemma5}
Let $p_1$, $p_2$ and $p_3$ be three points in a normed plane $\mathbb{M}^2$ such that $\|p_i-p_j\|\leq 1$ ($i,j\in\{1,2,3\}$) and
\begin{itemize}
\item there exists a point $x_{12}$ and a minimal circular arc $\widehat{p_1p_2}$ contained in $S(x_{12},1)$ such that $p_3$ is an interior point of $B(x_{12},1)$,
\item there exists a point $x_{23}$ and a minimal circular arc $\widehat{p_2p_3}$ contained in $S(x_{23},1)$ such that $p_1$ is an interior point of $B(x_{23},1)$,
\item $p_1$, $p_2$ and $p_3$ are not in a line.
\end{itemize}
Then $p_1\notin \operatorname{conv}(\widehat{p_2p_3}, \overline{p_2p_3})$ and $p_3\notin\operatorname{conv}(\widehat{p_1p_2}, \overline{p_1p_2})$.
\end{lem}

\smallskip

\begin{proof} Let us assume that $p_3\in \operatorname{conv}(\widehat{p_1p_2}, \overline{p_1p_2})$. Since $p_1$ is an interior point of $B(x_{23},1)$, then $S(x_{23},1)$ meets $\overline{p_1p_2}$ or $\widehat{p_1p_2}$ in a point different from $p_2$. In any case, the arc $\widehat{p_1p_2}$ is not contained in $B(x_{23},1)$, in contradiction to  Lemma \ref{4}. For similar reasons, $p_1\notin \operatorname{conv}(\widehat{p_2p_3}, \overline{p_2p_3})$.
\end{proof}

\bigskip

\begin{prop} \label{proposition1}
Let $K=\{p_1,p_2,\dots,p_n\}$ be a finite set in  a normed plane $\mathbb{M}^2$ having diameter 1. Then
$$\operatorname{bh}(K)=\bigcap_{i=1}^k B(x_i,1),$$
where $B(x_i,1)$, $i=1,2,...,k$, are balls which contain $K$, and their spheres contain some minimal arcs meeting points of $K$.
\end{prop}
%********************************************************\\
%I think that if $\|p_i-p_j\|\leq 1$, then through $p_i$ and $p_j$ there is a unique arc of radius 1. Why: if $C(x, 1)$ and $C(y, 1)$ are two different %circles through $p_i$ and $p_j$, then $C(x, 1)\cap C(y, 1)$ consists of two parallel to the line $\langle x, y\rangle$ segments (which can be %degenerate) being symmetric with respect to the midpoint of $[x, y]$. If the midpoint of $[x, y]$ coincides with the origin $o$, denote these two %segments by $[a, b]$ and $[-b, -a]$. The case when $p_i\in [a, b]$ and $[-b, -a]$ is impossible, because then after the translation though %$-\frac{1}{2}(p_i+p_j)$ we see that $\|p_i-p_j\|=2$, a contradiction. Thus we see that all the arcs of radius 1 with endpoints $p_i$ and $p_j$ %coincide with the line segment $[p_i, p_j]$.\\
%*********************************************************************

\smallskip

\begin{proof}

%It is clear that
%$$
%\operatorname{bh}(K)\subset \bigcap_{i=1}^k B(x_i,1).
%$$

%t y Lemma \ref{4}, any arc of radius 1 meeting points of $K$ is contained in any ball $B(x,1)\subset K$. Then
%$$\operatorname{conv}(\bigcup_{i,j=1}^n \widehat{p_ip_j}) \subseteq \operatorname{bh}(K).$$

%Let us consider a point $y\in \operatorname{bh}(K)$. Let $H$ be the set of all balls of radius 1 such that their boundary contains an arc meeting %points of $K$. Then
%$$y\in \operatorname{bh}(K)\ = \ \bigcap_{K\subset B(x,1)} B(x,1) \ \subset \
%\bigcap_{K\subset B(x,1)\subset H}B(x,1).
%$$
%If $y$ is from  an arc meeting two points of $K$, then $y\in \operatorname{conv}(\bigcup_{i,j=1}^n \widehat{p_ip_j}).$

%If $y$ is not from  an arc meeting two points of $K$, we choose a line $r$ which contains $y$. There exist two points $y_1$ and $y_2$ in $r$ such that %$y\in [y_1,y_2]$ and $y_i$ is an arc meeting points of $K$ for $i=1,2$.\\
%**********************************************************************\\
%Why $y_i$ are on an arc  meeting points of $K$ (see the attached figure)?  Or the line $r$ has a special position?\\
%****************************************************************************

% Since $y_i\in \operatorname{conv}(\bigcup_{i,j=1}^n \widehat{p_ip_j})$ for $i=1,2$, then $y\in \operatorname{conv}(\bigcup_{i,j=1}^n %\widehat{p_ip_j}).$

%\hfill$\square$

We fix the clockwise orientation of a closed curve in $\mathbb{M}^2$ as \emph{ negative orientation} of that curve.  Jung's Theorem (see, e.g., \cite{Bo-Ma})
guarantees the existence of a ball $B(x_1,1)$ such that $K\subset B(x_1,1)$. After translating  and renaming the points if  necessary, we may assume that
 \begin{itemize}
 \item $S(x_1,1)$ contains two points $p_1, p_2\in K$,
 \item the circular arc starting in $p_1$  with negative orientation and  ending in $p_2$ is a minimal circular arc,
 \item there is no other minimal circular arc in $S(x_1,1)$ meeting points of $K$ and being larger than the minimal circular arc meeting $p_1$ and $p_2$.
 \end{itemize}

The set $K$ has diameter 1, hence  $K$ is contained in $B(p_1,1)\cap B(p_2,1)\cap B(x_1,1)$. Starting in $z=x_1$, we move $z$ along $S(p_2,1)$ in the negative direction. Let $x_2$ denote  the first position of $z$ such that one of the following conditions is verified:
 \begin{enumerate}
 \item There is a new point $p_3\in S(x_2,1)$ such that
     \begin{itemize}
    \item the circular arc in $S(x_2,1)$ starting in $p_2$  with negative orientation and  ending in $p_3$ is a minimal circular arc,
    \item there is no other minimal circular arc in $S(x_2,1)$ meeting points of $K$ and being larger than the minimal circular arc meeting $p_2$ and $p_3$,
    \end{itemize}
 \item $p_1 \in S(x_2,1)$ with $x_2$ from  the other half-plane defined by the line $\langle p_1, p_2\rangle$.
 \end{enumerate}

 In both cases, we consider the set $A=B(x_1,1)\cap B(x_2,1)$. Since $z$ moves continuously in $S(p_2,1)$, $A$ contains $K$.

If $p_1\in S(x_2,1)$ with $x_2$ in  the other half-plane defined by  $\langle p_1, p_2\rangle$, then $A$ is the ball hull of $K$ by
$$\bigcap_{K\subset B(x,1)} B(x,1) \subset A,$$
since $A$ can be represented as intersection of balls $B(x_i,1)$ which contain $K$. Moreover, the boundary of $A$ is generated by minimal arcs meeting points of $K$. Thus $$A\subset  \bigcap_{K\subset B(x,1)} B(x,1).$$

If there exists a new point $p_3\in S(x_2,1)$, $K$ is contained in $B(p_1,1)\cap B(p_2,1)\cap B(p_3,1)\cap B(x_1,1)\cap B(x_2,1)$. Starting in $z=x_2$, we move $z$ along $S(p_3,1)$ in the negative direction. Let $x_3$ be the first value of $z$ such that one of the following conditions is verified:
 \begin{enumerate}
 \item There is a new point $p_4\in S(x_3,1)$ such that
 \begin{itemize}
    \item the circular arc in $S(x_3,1)$ starting in $p_3$ with negative orientation and  ending in $p_4$  is a minimal circular arc,
    \item there is no other minimal circular arc in $S(x_3,1)$ meeting points of $K$ and being larger than the minimal circular arc meeting $p_3$ and $p_4$.
    \end{itemize}
  \item $p_1 \in S(x_3,1)$.
\end{enumerate}
In both cases, we consider the set $A=B(x_1,1)\cap B(x_2,1)\cap B(x_3,1)$. Since $z$ moves continuously in $S(p_3,1)$, $A$ contains $K$.

If $p_1\in S(x_3,1)$, then $A$ is the ball hull of $K$ by
$$
\bigcap_{K\subset B(x,1)} B(x,1) \subset A,
$$
since $A$ is an intersection of balls $B(x_i,1)$ which contain $K$. Moreover, the boundary of $A$ is generated by minimal arcs meeting points of $K$. Thus $$A\subset  \bigcap_{K\subset B(x,1)} B(x,1).$$

If there is a new point $p_4$ in the above conditions, the process  continues  in the same way, and it is finite because the number of points $p_i$ is finite. In this process, starting in $p_j$, it  is not possible to get a point $p_i$ with $1<i<j$, due to Lemma \ref{lemma5}. At the end of the process, we obtain that the set
$A=\bigcap_{i=1}^kB(x_i,1)$
is the ball hull of $K$ because it is the intersection of unit balls which contain $K$, and its boundary is generated by minimal arcs meeting points of $K$.

\end{proof}

\bigskip

From Lemma \ref{4} and the proof of Proposition \ref{proposition1} we obtain

\begin{theo}
Let $K=\{p_1,p_2,\dots,p_n\}$ be a finite set in  a normed plane  $\mathbb{M}^2$ having  diameter 1. Let  $\widehat{p_ip_j}$ denote a minimal circular  arc of radius 1 meeting $p_i$ and $p_j$. Let $H$ be the set of all balls of radius 1 such that their boundary contains a circular  arc meeting points from $K$. Then
$$\operatorname{bh}(K)=\bigcap_{K\subset B(x,1)\subset H}B(x,1)=\operatorname{conv}(\bigcup_{i,j=1}^n \widehat{p_ip_j}).$$
\end{theo}

\vspace{1cm}
\begin{tabular}{l}
Pedro Mart\'{\i}n\\
Departamento de Matem\'{a}ticas,\\
Universidad de Extremadura,\\
06006 Badajoz, Spain\vspace{0.1cm}\\
E-mail: pjimenez@unex.es
\end{tabular}\vspace{0.3cm}

\begin{tabular}{l}
Horst Martini\\
Fakult\"at f\"ur Mathematik, TU Chemnitz\\
D-09107 Chemnitz, GERMANY\vspace{0.1cm}\\
E-mail: horst.martini@mathematik.tu-chemnitz.de
\end{tabular}\vspace{0.3cm}

\begin{tabular}{l}
Margarita Spirova\\
Fakult\"at f\"ur Mathematik, TU Chemnitz\\
D-09107 Chemnitz, GERMANY\vspace{0.1cm}\\
E-mail: margarita.spirova@mathematik.tu-chemnitz.de
\end{tabular}

\end{document}